\documentclass[12pt]{amsart}
\usepackage{amsfonts, amssymb, amsmath, amsthm, latexsym, array}
\usepackage{fullpage}
\usepackage{verbatim}

\input{xy}
\xyoption{all}

\newtheorem{thm}{Theorem}
\newtheorem{lm}{Lemma}

\newtheorem{prop}{Proposition}

\theoremstyle{remark}
\newtheorem{ex}{Example}

\theoremstyle{definition}

\renewcommand{\iff}{if and only if }

\newcommand{\gt}{\mathfrak}
\newcommand{\cp}{\mathbb C}

\newcommand{\GL}{{\rm GL}}

\newcommand{\ind}{{\rm ind\,}}

\newcommand{\rk}{{\rm rk\,}}
\newcommand{\ch}{{\rm char\,}}
\newcommand{\Lie}{{\rm Lie\,}}
\newcommand{\Ker}{{\rm Ker\,}}

\renewcommand{\Im}{{\rm Im\,}}
\newcommand{\Aut}{{\rm Aut}}
\newcommand{\Ann}{{\rm Ann}}
\newcommand{\ad}{{\rm ad }}

\begin{document}
\begin{center}
{\Large\bf The centralisers of nilpotent elements\\ in the classical Lie
algebras
}
\vskip1ex
{\sc O.S.~Yakimova}
\end{center}

\section*{Introduction}
Let
$\gt g$ be a Lie algebra over a field
$\mathbb K$. Consider the coadjoint representation
$\ad^*(\gt g)$.
{\it The index of} $\gt g$ is the minimum of dimensions
of stabilisers $\gt g_\alpha$ over all covectors
$\alpha\in\gt g^*$
$$
\ind\gt g=\min\limits_{\alpha\in\gt g^*}\dim\gt g_\alpha.
$$
The definition of index goes back to Dixmier
\cite[11.1.6]{Di}. This notion is important in
Representation Theory and also in Invariant Theory.
By Rosenlicht's theorem \cite{R}, generic
orbits of an arbitrary action of a linear algebraic group on
an irreducible algebraic variety are separated by rational
invariants; in particular, $\ind\gt g=\mbox{tr.deg\,}\mathbb
K(\gt g^*)^G$.

The index of a reductive algebra equals its rank. Computing
the index of an arbitrary Lie algebra seems to be a wild
problem. However, there is a number of interesting results
for several classes of non-reductive subalgebras of reductive
Lie algerbas. For instance, parabolic subalgebras and their
``relatives'' (nilpotent radicals, seaweeds) are considered
in \cite{al}, \cite{Dima1}, \cite{yu}. The centralisers of
elements form another interesting class of subalgebras. The last
topic is closely related to the theory of
integrable Hamiltonian systems.

Let $G$ be a semisimple Lie group (complex or real), $\gt
g=\Lie G$, and $Gx$ an orbit of a covector $x\in\gt g^*$.
Let $\gt g_x$ denote the stabiliser of $x$. It is well-known
that the orbit $Gx$ possesses a $G$-invariant symplectic
structure. There is a family of commuting with respect to a
Poisson bracket polynomial functions on $\gt g^*$
constructed by the argument shift method such that its
restriction to $Gx$ contains $\frac{1}{2}\dim(Gx)$
algebraically independent functions \iff $\ind\gt
g_x=\ind\gt g$, \cite{bol}.

\vskip1ex
\noindent {\bf Conjecture }{(\'Elashvili)}. {\it Let
$\gt g$ be a reductive Lie algebra. Then  $\ind\gt
g_x=\ind\gt g$ for each covector $x\in\gt g^*$.}

\vskip1ex
\noindent
Recall that if $\gt g$ is reductive,  then the
$\gt g$-modules $\gt g^*$ and $\gt g$ are isomorphic.
In particular, it is enough to prove
the ``index conjecture'' for stabilisers
of vectors $x\in\gt g$.

Given $x\in\gt g$, let  $x=x_s+x_n$ be the Jordan
decomposition. Then $\gt g_x=(\gt g_{x_s})_{x_n}$. The
subalgebra $\gt g_{x_s}$  is reductive and contains a Cartan
subalgebra of $\gt g$. Hence, $\ind\gt g_{x_s}=\ind\gt g=\rk
\gt g$. Thus, a verification of the "index conjecture" is
reduced to the computation of 
 $\ind\gt g_{x_n}$ for
nilpotent elements $x_n\in\gt g$. Clearly, we can restrict
ourselves to the case of simple $\gt g$.

Note that if
$x$ is a regular element, then the stabiliser $\gt g_x$ is commutative
and of dimension $\rk\gt g$.
The ``index conjecture'' was proved
for subregular nilpotents and nilpotents of height 2
\cite{Dima2}, and also for nilpotents of height
3  \cite{Dima3}. (The height of a
nilpotent element $e$ is the maximal number
$m$ such that $(\ad\, e)^m\ne 0$.)
Recently, \'Elashvili's conjecture was proved by Charbonnel
\cite{char} for $\mathbb K=\cp$.

In the present article, we prove in an elementary way,
that for any nilpotent element $e\in\gt g$ of a simple
classical Lie algebra the index of
$\gt g_e$ equals the rank of $\gt g$.
We assume that the ground field $\mathbb K$
contains at least $k$ elements, where $k$ is the number of
Jordan blocks of a nilpotent element $e\in\gt g$. For the
orthogonal and symplectic algebras, it is also assumed that
$\ch\mathbb K\ne 2$. Note that if a reductive Lie algebra
$\gt g$ does not contain exceptional ideals,
then $\gt g_{x_s}$ has the same property.
Thus, the ``index conjecture'' is proved for
the direct sums of classical algebras.

By Vinberg's inequality, which is  presented in
\cite[Sect.\,1]{Dima2}, we have $\ind\gt g_x\ge\ind\gt g$
for each element $x\in\gt g^*$.
It remains to prove the opposite 
inequality. 
To this end, it suffices to find $\alpha\in(\gt g_x)^*$ such
that the dimension of its stabiliser in $\gt g_x$ is at most
$\rk\gt g$. For $\gt g=\gt{gl}(V)$ and $\gt g=\gt{sp}(V)$,
we explicitly indicate such a point $\alpha\in\gt g_e^*$. In
case of the orthogonal algebra, the proof is partially based
on induction. 

In the last two sections, $\mathbb K$ is assumed to be
algebraically closed and of characteristic zero. It is shown
that the stabilisers $(\gt g_e)_\alpha$ constructed for $\gt
g=\gt{gl}(V)$ and $\gt g=\gt{sp}(V)$ are generic stabilisers
for the coadjoint representation of $\gt g_e$. For the
orthogonal case, we give an example of a nilpotent element
$e\in\gt{so}_8$ such that the coadjoint action of $\gt g_e$
has no generic stabiliser. Similar results for parabolic
and seaweed subalgebras of simple Lie algebras were obtained
by Panyushev and also by Tauvel and Yu. In \cite{yu} there
is an example of a parabolic subalgebra of $\gt{so}_8$
having no generic stabilisers for the coadjoint
representation. The affirmative answer for series $A$ and
$C$ is obtained by Panyushev in \cite{Dima4}.

In the last section, we consider the commuting variety of
$\gt{g}_e$ and its relationship with the commuting variety of
triples of matrices.

\par\medskip
This research was supported in part by CRDF grant
RM1--2543--MO--03.

\section{Preliminaries}
Suppose
$\gt g$ is a simple classical Lie algebra or
a general linear algebra.
Let $e\in\gt g$ be a nilpotent element and
$\gt z(e)$ its centraliser in $\gt g$.
Note that there is no essential difference
between $\gt g=\gt{gl}(V)$ and $\gt g=\gt{sl}(V)$.
However, the first case is more suitable for calculations.
In case of orthogonal and symplectic algebras,
we need some facts from the theory of symmetric spaces.

Let $(\ ,\ )_V$ be a non-degenerate
symmetric or skew-symmetric form on a finite dimensional
vector space $V$ given by a matrix $J$, i.e.,
$(v,w)_V=v^tJw$, where the symbol $\phantom{,}^t$ stands for
the transpose. The elements of $\gt{gl}(V)$ preserving
$(\phantom{,},\phantom{,})_V$ are exactly the fixed vectors
$\gt{gl}(V)^\sigma$ of the involution
$\sigma(\xi)=-J\xi^t J^{-1}$. There
is the $\gt{gl}(V)^\sigma$-invariant decomposition
$\gt{gl}(V)=\gt{gl}(V)^\sigma\oplus\gt g_1$.
The elements of
 $\gt g_1$ multiply the form
$(\phantom{,},\phantom{,})_V$ by $-1$, i.e.,
$(\xi v,w)_V=(v,\xi w)_V$ for every $v,w\in V$.

Set $\gt g=\gt{gl}(V)^\sigma$, and let $e\in\gt g$ be a nilpotent
element. Denote by $\gt z(e)$ and $\gt z_{\gt{gl}}(e)$
the centralisers of $e$ in
$\gt g$ and $\gt{gl}(V)$, respectively.
Since $\sigma(e)=e$,
$\sigma$ acts on $\gt z_{\gt{gl}}(e)$.
Clearly, $\gt z_{\gt{gl}}(e)^{\sigma}=\gt z(e)$.
This yields the
$\gt z(e)$-invariant decomposition
$\gt z_{\gt{gl}}(e)=\gt z(e)\oplus\gt z_{1}$.
Given $\alpha\in\gt z_{\gt{gl}}(e)^*$,
let $\tilde{\alpha}$ denote its restriction
to $\gt z(e)$.

\begin{prop}\label{1} Suppose $\alpha\in\gt z_{\gt{gl}}(e)^*$
and $\alpha(\gt z_{1})=0$. Then
$\gt z(e)_{\tilde{\alpha}}=\gt z_{\gt{gl}}(e)_\alpha\cap\gt z(e)$.
\end{prop}
\begin{proof} Take $\xi\in\gt z(e)$. Since
$[\xi, \gt z_{1}]\subset \gt z_{1}$, 
$\alpha([\xi, \gt z(e)])=0$
\iff $\alpha([\xi, \gt z_{\gt{gl}}(e)])=0$.
In particular, $\gt z(e)_{\tilde{\alpha}}=\gt z(e)_\alpha$.
\end{proof}

Suppose $\gt h$ is a Lie algebra and
$\tau\in\Aut\gt h$ an involution, which defines the decomposition
$\gt h=\gt h_0\oplus\gt h_1$. Each point
$\gamma\in\gt h_0^*$ determines a skew-symmetric
2-form $\hat\gamma$ on $\gt h_1$ by
$\hat\gamma(\xi,\eta)=\gamma([\xi,\eta])$.

\begin{lm} In the above notation, we have
$\ind\gt h\le\ind\gt h_0+\min\limits_{\gamma\in\gt h_0^*}\dim(\Ker\hat\gamma)$.
\end{lm}

\begin{proof} Consider $\gamma$ as a function on
$\gt h$, which is equal to zero on $\gt h_1$. Then $\gt
h_{\gamma}=(\gt h_0)_{\gamma}\oplus(\gt h_\gamma\cap\gt
h_{1}) =(\gt h_0)_{\gamma}\oplus(\Ker\hat\gamma)$. We have
$\dim(\gt h_0)_\gamma=\ind\gt h_0$ for generic points (=\,
points of some Zariski open subset $U_1\subset\gt h_0^*$).
The points of $\gt h_0^*$, where $\Ker\hat\gamma$ has the
minimal possible dimension, form another open subset, say
$U_2\subset\gt h_0^*$. For the points of the intersection
$U_1\cap U_2$, the dimension of the stabiliser in $\gt h$
equals the required sum.
\end{proof}


\section{General linear algebra}

Consider a nilpotent element
$e\in\gt{gl}(V)$, where
$V$ is an $n$-dimensional vector space over
$\mathbb K$. Denote by $\gt z(e)$ the centraliser of $e$.
Let us show that the index of $\gt z(e)$ equals $n$.

Let $k$ be a number of Jordan blocks of $e$ and $W\subset V$
a $k$-dimensional complement of $\Im e$ in $V$. Denote by
$d_i+1$ the dimension of $i$-th Jordan block. Choose a basis
$w_1, w_2, \ldots, w_k$ in $W$, where $w_i$ is a generator
of an $i$-th Jordan block, i.e., the vectors $e^{s_i}w_i$
with $1\le i\le k$, $0\le s_i\le d_i$ form a basis of $V$.
Let $\varphi\in\gt z(e)$. Since
$\varphi(e^sw_i)=e^s\varphi(w_i)$, the map $\varphi$ is
completely determined by its values on $w_i$,
$i=1,\ldots,k$. Each value $\varphi(w_i)$ can be written as
$$
\varphi(w_i)=\sum_{j,s} c_i^{j,s}e^s w_j, \mbox{ where }
c_i^{j,s}\in\mathbb K. \eqno (1)
$$
That is,
$\varphi$ is completely determined by the coefficients
$c_i^{j,s}=c_i^{j,s}(\varphi)$. Note that
$\varphi\in\gt z(e)$ preserves the space of each Jordan block
\iff
$c_i^{j,s}(\varphi)=0$ for $i\ne j$.

The centraliser $\gt z(e)$ has a basis
$\{\xi_i^{j,s}\}$, where
$$
\begin{array}{ll}
\left\{
\begin{array}{l}
\xi_i^{j,s}(w_i)=e^sw_j, \\
\xi_i^{j,s}(w_t)=0 \enskip \mbox{for } t\ne i, \\
\end{array}\right. &
\left[
\begin{array}{l}
d_j-d_i\le s\le d_j \enskip \mbox{for } d_j\ge d_i, \\
0\le s\le d_i \enskip \mbox{for } d_j \end{array}\right. \\
\end{array}
$$

Consider a point
$\alpha\in\gt z(e)^*$ defined by the formula
$$
\alpha(\varphi)=\sum\limits_{i=1}^k a_i\cdot c_i^{i,d_i},
\enskip a_i\in\mathbb K,
$$
where $c_i^{j,s}$ are the coefficients
of $\varphi\in\gt z(e)$ and $\{a_i\}$ are non-zero pairwise distinct
numbers. We have
$\alpha(\xi_i^{j,s})=a_i$ if $i=j$, $s=d_i$ and zero otherwise.

\begin{thm} The stabiliser
$\gt z(e)_\alpha$ of $\alpha$ in $\gt z(e)$ consist of all
maps preserving the Jordan blocks, i.e., $\gt z(e)_\alpha$
is the linear span of the vectors $\xi_i^{i,s}$.
\end{thm}

\begin{proof}
Suppose $\varphi\in\gt z(e)$ is defined by formula $(1)$.
(Some of $c_i^{j,s}$ have to be zeros, but this is immaterial here).
For each basis vector $\xi_i^{j,b}$, we have
$$
\alpha([\varphi, \xi_i^{j,b}])=\alpha(\sum\limits_{t,s}c_t^{i,s}\xi_t^{j,s+b}-
\sum\limits_{t,s}c_j^{t,s}\xi_i^{t,s+b})=
a_j\cdot c_j^{i,d_j-b}-
a_i\cdot c_j^{i,d_i-b} \ .
$$
The element
$\varphi$ lies in  $\gt z(e)_\alpha$ \iff
$\alpha([\varphi, \xi_i^{j,b}])=0$ for all
$ \xi_i^{j,b}$.

Note that if $\varphi$ preserves the Jordan blocks, i.e.,
$c_i^{j,s}=0$ for $i\ne j$, then $\alpha([\varphi, \gt
z(e)]=0$. Let us show that $\gt z(e)_\alpha$
contains no other elements.
Assume that $c_i^{j,s}\ne 0$ for some
$i\ne j$. We have three different possibilities:
$d_i<d_j$, $d_i=d_j$ and $d_i>d_j$.

If $d_j\le d_i$, then put $\xi(w_j)=e^{(d_i-s)}w_i$ and
$\xi(w_t)=0$ for $t\ne j$. It should be noted that
$0\le s\le d_j\le d_i$, hence the expression $e^{(d_i-s)}$
is well defined. One has to check that
$e^{d_j+1}(\xi(w_j))=0$. Adding the
powers of $e$, we get
$e^{d_j+1}(\xi(w_j))=e^{d_j+1+d_i-s}w_i=e^{d_j-s}(e^{d_i+1}w_i)=0$.
We have $\alpha([\varphi,\xi])=a_i\cdot c_i^{j,s}-a_j\cdot
c_i^{j, d_j-d_i+s}$. In case $d_j=d_i$, we obtain
$(a_i-a_j)\cdot c_i^{j,s}\ne 0$. If $d_j>d_i$, then
$s>d_j-d_i+s$. Choose the minimal $s$ such that
$c_i^{j,s}\ne 0$. For this choice, we get
$\alpha([\varphi,\xi])=a_i\cdot c_i^{j,s}\ne 0$.

Suppose now that $d_j>d_i$ and  $s$ is the minimal number
such that $c_i^{j,s}\ne 0$. Set $\xi(w_j)=e^{(d_j-s)w_i}$ and
$\xi(w_t)=0$ for $t\ne j$. As in the previous case, we have
$0\le s\le d_j$. In particular, $d_j-s\ge 0$, $(d_j+1+d_j-s)>d_i+1$
and, thereby, $e^{d_j+1}(\xi(w_j))=0$. We obtain
$$
\alpha([\xi,\varphi])=a_j\cdot c_i^{j,s}-a_i\cdot
c_i^{j, d_i-d_j+s}=a_j\cdot c_i^{j,s}\ne 0.
$$
Here $c_i^{j, d_i-d_j+s}=0$, since $d_i-d_j+s<s$.
\end{proof}

\vskip0.2ex
\noindent
{\bf Corollary.} {\it The index of $\gt z(e)$ equals $n$.}

\begin{proof} The stabiliser
$\gt z(e)_\alpha$ consist of all maps preserving Jordan
blocks. In particular, it has dimension $n$. Hence, $\ind\gt
z(e)\le n$.
On the other hand, it follows from
Vinberg's inequality that $\ind\gt z(e)\ge n=\rk\gt{gl}(V)$.
\end{proof}

Let us give another proof of the inequality
$\ind\gt z(e)\le n$. 

\begin{ex}
Let $e\in\gt{gl}_n$ be a nilpotent element and $\gt h=\gt
z(e)$ the centraliser of $e$. We may assume that the first
Jordan block of $e$ is of maximal dimension. Then
$V=V_{d_1+1}\oplus V_{\mbox{\footnotesize oth}}$ and
$e=e_1+e_2$, where $V_{d_1+1}$ is the space of the first
Jordan block and $V_{\mbox{\footnotesize oth}}$ is the space
of all other Jordan blocks; $e_1\in\gt{gl}_{d_1+1}$,
$e_2\in\gt{gl}_{n-d_1-1}$.
Let $\tau\in \gt{gl}(V)$ be the conjugation by a diagonal
matrix of order two such that
$\gt{gl}(V)^\tau=\gt{gl}_{d_1+1}\oplus\gt{gl}_{n-d_1-1}$.
The involution $\tau$ acts on  $\gt h=\gt z(e)$ and induces
the decomposition $\gt h=\gt h_0\oplus\gt h_{1}$, where $\gt
h_0=\gt z(e_1)\oplus\gt z(e_2)$ (the centralisers are
considered in the algebras $\gt{gl}_{d_1+1}$ and
$\gt{gl}_{n-d_1-1}$, respectively). Assume that the ``index
conjecture'' is true for all $m<n$; in particular, $\ind\gt
z(e_2)=n-d_1-1$. The subalgebra $\gt z(e_1)$ is commutative
and its index equals $d_1+1$. According to Lemma~1, $\ind\gt
z(e)\le\ind(\gt z(e_1)\oplus\gt z(e_2))+
\min\limits_{\gamma\in\gt h_0^*}\dim(\Ker\hat\gamma)\le
n+\min\limits_{\gamma\in\gt z(e_1)^*}\dim(\Ker\hat\gamma)$.
Now, we make a special choice for $\gamma$. Set
$\gamma(\xi_1^{1,d_1})=1$ and $\gamma(\xi_i^{j,s})=0$ for
all other $\xi_i^{j,s}$. The subspace $\gt h_{1}$ is
generated by the vectors $\xi_i^{1,s}$ and $\xi_1^{i,s}$
with $i\ne 1$. We have
$$
\left\{
\begin{array}{l}
\hat\gamma(\xi_1^{i,s},\xi_i^{1,d_1-s})=1, \\
\hat\gamma(\xi_1^{i,s},\xi_i^{1,b})=0 \mbox{ if } s+b\ne d_1, \\
\hat\gamma(\xi_1^{i,s},\xi_1^{i,b})=
\hat\gamma(\xi_i^{1,s},\xi_i^{1,b})=0.\\
\end{array}
\right.
$$
The form $\hat\gamma$ defines a non-degenerate
pairing between the 
spaces
$U_1^i:=\left<\xi_1^{i,s}|0\le s\le d_i\right>$
and $U_i^1:=\left<\xi_i^{1,s}|d_1-d_i\le s\le d_1\right>$.
Hence, $\hat\gamma$ is non-degenerate and
$\ind\gt z(e)\le n$.
\end{ex}

\section{Symplectic algebra}

In this section $\gt g=\gt{sp}_{2n}=\gt{sp}(V)$, where $V$
in an $2n$-dimensional vector space over $\mathbb
K$. As above, $e\in\gt{sp}_{2n}$ is a nilpotent element and
$\gt z(e)\subset\gt g$ is the centraliser of $e$. 
Let $\{w_i\}$ be generators of Jordan blocks
associated with $e$. We may assume that  
the space of each 
even-dimensional Jordan block is orthogonal to the space of 
all other Jordan blocks. If $d_i$ is even, then the restriction 
of $\gt{sp}_{2n}$-invariant form $(\,,\,)_V$ on the space of 
$i$-th Jordan block is zero. One can choose generators $\{w_i\}$ 
such that the odd-dimensional blocks are partitioned in pairs $(i,i')$, 
where $i'$ is the number of the unique Jordan block which is not
orthogonal to the $i$-th one. Note that $d_{i'}=d_i$.

Let $\gt z_{\gt{gl}}(e)$ be the centraliser of $e$ in
$\gt{gl}_{2n}$. Recall that $\gt z(e)=\gt
z_{\gt{gl}}(e)^\sigma\oplus\gt z_1$, where $\sigma$ is an
involutive automorphism of $\gt{gl}_{2n}$. For elements of
$\gt z_{\gt{gl}}(e)$ we use notation introduced in the
previous section.

Let $\alpha\in\gt
z_{\gt{gl}}(e)^*$ be a function determined just like in the
previous case:
$$
\alpha(\varphi)=a_1\cdot
c_1^{1,d_1}+a_2\cdot c_2^{2,d_2}+\ldots+a_{2n}\cdot c_{k}^{k,d_{k}},
$$
where $\varphi$ is given by its coefficients
$c_i^{j,s}$,
and $\{a_i\}$ are pairwise distinct non-zero numbers
with $a_{i'}=-a_i$.

\begin{lm} In the above notation, we have $\alpha(\gt z_{1})=0$.
\end{lm}

\begin{proof}
Assume that there is
$\psi\in\gt z_{1}$ such that $\alpha(\psi)\ne 0$.
Then there is a non-zero coefficient $c_i^{i,d_i}$ of $\psi$.
Recall that $\sigma(\psi)=-\psi$. The element $\psi$
multiplies the $\gt{sp}_{2n}$-invariant skew-symmetric form
$(\phantom{,},\phantom{,})_V$ by $-1$, in particular,
$(\psi(w_i),v)_V=(w_i,\psi(v))_V$ for each vector $v\in V$.
Clearly, $\psi(w_i)$ and $w_i$ have to be orthogonal with
respect to the skew-symmetric form. If $d_i$ is odd, then
$(w_i,e^{d_i}w_i)_V\ne 0$, hence, $c_i^{i,d_i}=0$. If on the
contrary $d_i$ is even, then
$$
\begin{array}{l}
c_i^{i,d_i}(e^{d_i}w_i, w_{i'})_V=(\psi(w_i), w_{i'})_V=
(w_i, \psi(w_{i'}))_V=c_{i'}^{i',d_i}(w_i,e^{d_i}w_{i'})_V= \\
\qquad \qquad \qquad\qquad \qquad\qquad
=(-1)^{d_i}c_{i'}^{i',d_i}(e^{d_i}w_i, w_{i'})_V=
c_{i'}^{i',d_i}(e^{d_i}w_i, w_{i'})_V. \\
\end{array}
$$
Hence, $c_i^{i,d_i}=c_{i'}^{i',d_i}$. Combining this
equality with defining formula of $\alpha$ we get a sum over
pairs of odd-dimensional blocks
$$
\alpha(\psi)=\sum\limits_{(i,i')}(a_i+a_{i'})c_i^{i,d_i},
$$
which is zero since $a_i=-a_{i'}$.
\end{proof}

Denote by
$\tilde\alpha$ the restriction of $\alpha$
to $\gt z(e)$.

\begin{thm} The dimension of the
stabiliser
$\gt z(e)_{\tilde\alpha}=\gt
z_{\gt{gl}}(e)_\alpha\cap\gt{sp}_{2n}$ equals $n$.
\end{thm}

\begin{proof}
The stabiliser of $\alpha$ in $\gt z_{\gt{gl}}(e)$ consist
of all maps preserving the spaces of the Jordan blocks. By
Proposition 1, $\gt z(e)_{\tilde{\alpha}}= \gt
z_{\gt{gl}}(e)_\alpha\cap\gt z(e)$. Describe the
intersection of $\gt z_{\gt{gl}}(e)_\alpha$ with the
symplectic subalgebra. If $w_i$ is a generator of an
even-dimensional block, then $\xi_i^{i,s}$ multiply the
skew-symmetric form by $(-1)^{s+1}$, i.e.,
$(\xi_i^{i,s}(e^bw_i),
e^tw_i)=(-1)^s(e^bw_i,\xi_i^{i,s}(e^tw_i))$. Consider a
space of a pair $(i,i')$ of odd-dimensional blocks. Set
$d:=d_i=d_{i'}$. Recall that $(w_i,
e^dw_{i'})=(-1)^s(e^sw_i, e^{d-s}w_{i'})=-(w_{i'},e^dw_i)$.
Since $(e^sw_i, e^{d-s}w_{i'})=(-1)^s(w_i, e^dw_{i'})$, the
elements $\xi_i^{i,s}+(-1)^{s+1}\xi_{i'}^{i',s}$ preserve
the skew-symmetric form, and the elements
$\xi_i^{i,s}+(-1)^{s}\xi_{i'}^{i',s}$ multiply it by $-1$.
From each even-dimensional block $i$ we get $(d_i+1)/2$
vectors, and from a pair $(i,i')$ we get $d_i+1$ vectors.
Thus the stabiliser of $\alpha$ in the whole of
$\gt{sp}_{2n}$ is an $n$-dimensional subalgebra.
\end{proof}

\section{The orthogonal case}

In this section $\gt g=\gt{so}_n$. As above
$e\in\gt{so}_{n}$ is a nilpotent element, $\gt z(e)$ is the
centraliser of $e$ in $\gt g$. 
Let $\{w_i\}$ be generators of Jordan blocks
associated with $e$. We may assume that  
the space of each 
odd-dimensional Jordan block is orthogonal to the space of 
all other Jordan blocks. If $d_i$ is odd, then the restriction 
of $\gt{so}_{n}$-invariant form $(\,,\,)_V$ on the space of 
$i$-th Jordan block is zero. One can choose generators $\{w_i\}$ 
such that the even-dimensional blocks are partitioned in pairs $(i,i^*)$, 
where $i^*$ is the number of the unique Jordan block which is not
orthogonal to the $i$-th one. Note that $d_{i^*}=d_i$.

Like the symplectic algebra, the orthogonal algebra is a
symmetric subalgebra of $\gt{gl}_{n}$. Denote by $\sigma$
the involution defining it. Since $\sigma(e)=e$, we have
 $\gt
z_{\gt{gl}}(e)=\gt z(e)\oplus\gt z_1$
similarly to the symplectic case. 
If $d_i$ is even, set $i^*=i$.
Assume that
$(w_{i*},e^{d_i}w_i)_V=\pm 1$ and $(w_i,e^{d_i}w_i)_V=1$ for
$i=i^*$. The algebra $\gt z(e)$ is generated (as a vector
space) by the vectors
$\xi_i^{j,d_j-s}+\varepsilon(i,j,s)\xi_{j^*}^{i^*,d_i-s}$,
where $\varepsilon(i,j,s)=\pm 1$ depending on $i,j$ and $s$.
In its turn, the subspace $\gt z_{1}$ is generated by the
vectors
$\xi_i^{j,d_j-s}-\varepsilon(i,j,s)\xi_{j^*}^{i^*,d_i-s}$.
Recall that $(e^sw_i,e^{d_i-s}w_{i^*})_V\ne 0$ if $e^sw_i\ne
0$.

We give some simple examples of linear functions with zero
restrictions to $\gt z_1$. Let $\varphi\in\gt
z_{\gt{gl}}(e)$ be a linear map defined by Formula $(1)$.
Set $\beta_i(\varphi)=c_i^{i,d_i-1}$,
$\gamma_{i,j}(\varphi)=c_i^{j,d_j}$.

\begin{lm}\label{fso} If $i=i^*$, $j=j^*$,
$t\ne t^*$, then
functions $\beta_i$, $\gamma_{i,j}-\gamma_{j,i}$
and $\gamma_{t,t}+\gamma_{t^*,t^*}$ are
equal to zero on $\gt z_1$.
\end{lm}

\begin{proof} Suppose  $\psi\in\gt z_1$ is defined
by Formula $(1)$. Since $\sigma(\psi)=-\psi$ and
$(\psi(w_i),ew_i)_V=c_i^{i,d_i-1}(e^{d_1-1}w_i,ew_i)_V$,
 we have
$$
(\psi(w_i),ew_i)_V=(w_i,\psi(ew_i))_V=
(w_i,e\psi(w_i))_V=-(ew_i,\psi(w_i))_V=
-c_i^{i,d_i-1}(ew_i,e^{d_1-1}w_i)_V.
$$
The form $(\phantom{,},\phantom{,})_V$ is symmetric and
$(ew_i,e^{d_i-1}w_{i})_V\ne 0$, hence
$\beta_i(\psi)=c_i^{i,d_i-1}=0$.

Similarly,
$$
\begin{array}{l}
c_i^{j,d_j}(e^{d_j}w_j,w_j)_V=
(\psi(w_i),w_j)_V=(w_i,\psi(w_j))_V=c_j^{i,d_i}(w_i,e^{d_i}w_i)_V;\\
c_t^{d_t,t}(e^{d_t}w_t,w_{t*})_V=(\psi(w_t),w_{t^*})_V=
(w_t,\psi(w_{t^*}))_V=c_{t^*}^{t^*,d_t}(w_t,e^{d_t}w_{t^*})_V.\\
\end{array}
$$
Recall that by our choice
$(e^{d_j}w_j,w_j)_V=(w_i,e^{d_i}w_i)_V=1$,
$(e^{d_t}w_t,w_{t*})_V=-(w_t,e^{d_t}w_{t^*})_V$. Hence
$c_i^{j,d_j}=c_j^{i,d_i}$, $c_t^{t,d_t}=-c_{t^*}^{t^*,d_t}$.
\end{proof}

Let us prove the inequality
$\ind\gt z(e)\le\rk\gt{so}_n$ by the
induction on 
$n$.
In the following two cases, the induction argument
does not go through. Therefore we consider them separately.

\noindent
{\it The first case.}
If $e\in\gt{so}_{2m+1}$ is a regular nilpotent
element, then
$\gt z(e)$ is a commutative $m$-dimensional algebra.

\noindent
{\it The second case.}
Let $e\in\gt{so}_{4d}$ be a nilpotent element
with two Jordan blocks of size $2d$ each.
Set  $\alpha(\varphi)=c_1^{1,2d-2}-c_2^{2,2d-2}$,
where $\varphi$ is defined by Formula $(1)$.
One can easily check that $\gt z(e)_\alpha$
has a basis
$\xi_1^{1,s}+(-1)^{s+1}\xi_2^{2,s}$ with $0\le s\le 2d-1$
and that $\dim\gt z(e)_\alpha=2d$.

Order the Jordan blocks of $e$ according to their dimensions
$d_1\ge d_2\ge\ldots\ge d_k$. Here
$d_i+1$ stands for the dimension of the
$i$-th Jordan block, similarly to the
case of $\gt{gl}_n$. Note that
the numbers
$n$ and $k$ have the same parity.
Assume that $k>1$ and if
$k=2$, then both Jordan blocks are odd dimensional.
Then we have the following three possibilities:

{\bf (1)} for some even number $2p<k$ the restriction of
$(\phantom{,},\phantom{,})_V$ to the space of the first $2p$
Jordan blocks is non-degenerate;

{\bf (2)} the number $d_i$ is even for $i=1,k$ and odd for
all other $i$;

{\bf (3)} the number $d_i$ is even \iff $i=1$.

\noindent 
Each of these three possibilities is considered
separately. In the first two cases we make an induction
step. In the third one a point $\alpha\in\gt z(e)^*$ is
given such that  $\dim\gt z(e)_\alpha\le\rk\gt{so}_n$.

\vskip0.1ex
 \noindent
  {\bf (1)} Suppose the space $V_{2m}$
of the first $2p$ Jordan blocks has dimension $2m$ and the
restriction of $(\phantom{,},\phantom{,})_V$ to $V_{2m}$ is
non-degenerate. Then $V=V_{2m}\oplus V_{\mbox{\small oth}}$,
$e=e_1+e_2$, where $e_1\in\gt{so}_{2m}$,
$e_2\in\gt{so}_{n-2m}$. Let $\tau$ be an involution of
$\gt{gl}_n$ corresponding to these direct sum, i.e.,
$\gt{gl}_n^\tau=\gt{gl}(V_{2m})
\oplus\gt{gl}(V_{\mbox{\small oth}})$. Set $\gt h=\gt z(e)$,
$\gt h_0=\gt z(e)^\tau$. Then $\gt h_0=\gt z(e_1)\oplus\gt
z(e_2)$, where the centralisers of $e_1$ and $e_2$ are taken
in $\gt{so}_{2m}$ and $\gt{so}_{n-2m}$, respectively.
By the inductive hypothesis, $\ind\gt z(e_1)=m$, $\ind\gt
z(e_2)=[n/2]-m$. Hence, $\ind\gt z(e)\le
[n/2]+\min\limits_{\gamma\in\gt h_0^*}
\dim(\Ker\hat\gamma)$. To conclude we have to point out a
function $\gamma\in\gt h_0^*$ such that $\hat\gamma$ is
non-degenerate. Recall that the involutions $\sigma$ and
$\tau$ commute with each other, preserve $e$ and determine the
decomposition $\gt z_{\gt{gl}}(e)=(\gt z(e_1)\oplus\gt
z(e_2)\oplus \gt h_1)\oplus\gt z_1$. If $\gamma(\gt
z_1)=\gamma(\gt h_1)=0$, then $\Ker\hat\gamma=(\gt
h_1\cap\gt z_{\gt{gl}}(e)_\gamma)$.

Divide odd-dimensional Jordan blocks into pairs  $(i,i')$
(it is assumed that $i,i'\le 2p$). Define a point $\gamma$
by
$$
\gamma(\varphi)=\sum\limits_{(i,i'),\, i,i'\le 2p}
(c_i^{i',d_{i'}}-c_{i'}^{i,d_i})+ \sum\limits_{j\le 2p,\,
(d_j+1) \mbox{{\small{ is even}}}} c_j^{j,d_j},
$$
where $\varphi\in\gt z_{\gt{gl}}(e)$ is given by its
coefficients  $c_i^{j,s}$.
The first summand is a sum of $(\gamma_{i,i'}-
\gamma_{i',i})$ over pairs of
odd-dimensional blocks, the second is the sum
of $(\gamma_{j,j}+\gamma_{j^*,j^*})$ over
pairs of even-dimensional blocks.
According to Lemma~\ref{fso}, both summands
are identical zeros on
$\gt z_1$. Moreover, by the definition
$\gamma(\gt h_1)=0$.

Set $j':=j$ for even-dimensional blocks. Assume that an
element $\psi\in\gt h_1$ determined by $(1)$  lies in the
kernel of $\hat\gamma$, i.e., $\gamma([\psi,\gt h_1])=0$.
Then $\gamma([\psi, \gt h])=\gamma([\psi,\gt
z_{\gt{gl}}(e)])=0$. Since $\psi\in\gt{so}_{n}$ and $\psi\ne
0$, we may assume that  $c_i^{j,s}\ne 0$ for some $j>2p\ge
i$. We have
$$
\gamma([\psi, \eta_i^{j',d_{j'}}])=\pm c_i^{j,s}\ne 0.
$$
Thus we have proved that
$\hat\gamma$ is non-degenerate and
$\ind\gt z(e)\le [n/2]$.

\vskip0.2ex 
\noindent 
{\bf (2)} Consider a decomposition
$V=V_{\mbox{\small oth}}\oplus V_{d_k+1}$, where the second
summand is the space of the smallest (odd-dimensional)
Jordan block and the first one is the space of all other 
blocks. As above $e=e_1+e_2$, where nilpotent element $e_2$
corresponds to the smallest (odd-dimensional) Jordan block.
We define an involution $\tau$, algebras $\gt z(e_1)$, $\gt
z(e_2)$, $\gt h_0$ and a subspace $\gt h_1$ in the same way
as in case {\bf (1)}. By the inductive hypothesis $\ind\gt
z(e_1)=[(n-d_k-1)/2]$, $\ind\gt z(e_2)=d_k/2$. Hence
$\ind\gt h_0=n/2-1$. By Lemma~1, $\ind\gt
z(e)\le\ind\gt h_0 + \min\limits_{\gamma\in\gt
h_0^*}\dim(\Ker\hat\gamma)$. Let $\gamma$ be the following
function
$$
\gamma(\varphi)=c_1^{1,d_1-1}+\sum\limits_{i=2}^{k-1} c_i^{i,d_i},
$$
where $\varphi$ is given by formula  $(1)$.
The first summand is $\beta_1$, the second summand
is a sum of $(\gamma_{j,j}+\gamma_{j^*,j^*})$
over pairs of even-dimensional blocks.
Due to Lemma~\ref{fso},
$\gamma(\gt z_1)=0$.
Suppose $\psi\in\gt h_1$ is given by its coefficients
$c_i^{j,s}$. Then
$$
\left\{
\begin{array}{l}
\gamma([\xi_1^{k,b},\psi])=c_k^{1,d_1-1-b},\\
\gamma([\xi_i^{k,b},\psi])=c_k^{i,d_i-b} \quad \mbox{ for } 1<i<k.\\
\end{array}\right.
$$
One can see that the kernel of
$\hat\gamma$ is one-dimensional and generated by
$(\xi_k^{1,d_1}-\xi_1^{k,d_k})$.
Hence,
$\ind\gt z(e)\le [n/2]-1+1=[n/2]$.

\vskip0.1ex
\noindent
{\bf (3)} In this case $n$ and $k$ are odd,
and
$e$ has a unique odd-dimensional Jordan block whose size is
maximal.
Assume that $k=2m+1$.
Enumerate the Jordan blocks by integers ranging from
$-m$ to $m$. Let the unique odd-dimensional
block has number zero. Suppose that pairs of blocks
$(-i,i)$ and $(-j,j)$ are orthogonal to each other if
$i\ne\pm j$, and dimensions of Jordan
blocks are increasing from  $-m$ to $0$ and decreasing from
$0$ to $m$, i.e., if $|i|\le |j|$, then $d_i\ge d_j$. Note 
that  $d_i=d_{-i}$. Such enumeration is shown on Picture 1.
Choose the generators $w_i$ of Jordan blocks such that
$i(w_i,e^{d_i}w_{-i})_V=|i|$ for
$i\ne 0$ and $(w_0,e^{d_0}w_0)_V=1$.

\begin{figure}[htb]
\setlength{\unitlength}{0.017in}
\begin{center}
\begin{picture}(150,110)(0,0)
\put(0,0){\line(1,0){150}}
\multiput(0,0)(150,0){2}{\line(0,1){20}}
\put(0,20){\line(1,0){10}}
\put(150,20){\line(-1,0){10}}
\multiput(10,20)(130,0){2}{\line(0,1){20}}
\multiput(30,40)(90,0){2}{\line(0,1){20}}
\multiput(40,60)(70,0){2}{\line(0,1){20}}
\multiput(60,80)(30,0){2}{\line(0,1){20}}
\multiput(70,100)(10,0){2}{\line(0,1){10}}
\multiput(10,40)(110,0){2}{\line(1,0){20}}
\multiput(30,60)(80,0){2}{\line(1,0){10}}
\multiput(40,80)(50,0){2}{\line(1,0){20}}
\multiput(60,100)(20,0){2}{\line(1,0){10}}
\put(70,110){\line(1,0){10}}

\qbezier[12](10,0),(10,10),(10,20) 
\qbezier[24](20,0),(20,20),(20,40)
\qbezier[24](30,0),(30,20),(30,40)
\qbezier[36](40,0),(40,30),(40,60)
\qbezier[48](50,0),(50,40),(50,80)
\qbezier[48](60,0),(60,40),(60,80)
\qbezier[60](70,0),(70,50),(70,100)
\qbezier[60](80,0),(80,50),(80,100)
\qbezier[48](90,0),(90,40),(90,80)
\qbezier[48](100,0),(100,40),(100,80)
\qbezier[36](110,0),(110,30),(110,60)
\qbezier[24](120,0),(120,20),(120,40)
\qbezier[24](130,0),(130,20),(130,40)
\qbezier[12](140,0),(140,10),(140,20)


\put(-1.5,-7){$\scriptstyle -m$}
\put(14,-7){$\ldots$}
\put(29,-7){$\scriptstyle -4$}
\put(39,-7){$\scriptstyle -3$}
\put(49,-7){$\scriptstyle -2$}
\put(59,-7){$\scriptstyle -1$}
\put(73.5,-7){$\scriptstyle 0$}
\put(83.5,-7){$\scriptstyle 1$}
\put(93.5,-7){$\scriptstyle 2$}
\put(103.5,-7){$\scriptstyle 3$}
\put(111.9,-7){$\ldots$}
\put(127.9,-7){$\scriptstyle m-1$}
\put(143.8,-7){$\scriptstyle m$}
\end{picture}
\end{center}
\caption{}\label{pikcha_A}
\end{figure}

Suppose $\varphi\in\gt z_{\gt{gl}}(e)$ is given by
Formula $(1)$. Consider the following point
$\alpha\in\gt z_{\gt{gl}}(e)^*$:
$$
\alpha(\varphi)=\sum\limits_{i=-m+1}^m c_{i-1}^{i,d_i}.
$$
One can check by direct computation that
$c_{i-1}^{i,d_i}(\psi)=-c_{-i}^{1-i,d_{1-i}}(\psi)$ for each
$\psi\in\gt z_1$ and, hence, $\alpha(\gt z_1)=0$. Let
$\tilde\alpha\in\gt z(e)^*$ be the restriction of $\alpha$.
Let us describe the stabiliser $\gt z(e)_{\tilde\alpha}=\gt
z_{\gt{gl}}(e)_\alpha\cap\gt z(e)$. Note that
$\alpha([\varphi,\xi_i^{j,s}])=c_{j-1}^{i,d_j-s}(\varphi)
-c_j^{i+1,d_{i+1}-s}(\varphi)$.

\begin{lm} Suppose  $\varphi\in\gt z(e)$ and
$\ad^*(\varphi)\alpha=0$. Then
$c_i^{j,s}=c_i^{j,s}(\varphi)=0$ for $i<j$.
\end{lm}

\begin{proof} Assume that the statement is wrong
and take a maximal $i$ for which there are
$j>i$ and $s$ such that $c_i^{j,s}\ne 0$. Because
$\varphi$ preserves $(\phantom{,},\phantom{,})_V$,
$c_{-j}^{-i, d_i-d_j+s}=\pm c_i^{j,s}\ne 0$.
Hence,      $-j\le i<j$, $j>0$, $|i|\le j$ and
$d_i\ge d_j$. Moreover,  $-j<(i+1)\le j$ and
$d_{i+1}\ge d_j$. Evidently,
$d_{i+1}-s\ge d_j-s\ge 0$ and there is an element
$\xi_j^{i+1,d_{i+1}-s}\in\gt z_{\gt{gl}}(e)$.
We have
$$
0=\alpha([\varphi,\xi_j^{i+1,d_{i+1}-s}])=
c_i^{j,s}-c_{i+1}^{j+1,\delta}=c_i^{j,s}.
$$
Here we do not give a precise value of $\delta$. Anyway all
coefficients $c_{i+1}^{j+1,b}$ are zeros, because
$j+1>i+1>i$. We get a contradiction. Thus the lemma is
proved.
\end{proof}

Let us say that $\varphi\in\gt z_{\gt{gl}}(e)$ has a {\it
step} $l$ whenever $c_i^{j,s}(\varphi)=0$  for $j\ne i+l$.
Each vector $\varphi\in\gt z_{\gt{gl}}(e)$ can be
represented as a sum
$\varphi=\varphi_{-2m}+\varphi_{-2m+1}+\ldots+\varphi_{2m-1}
+\varphi_{2m}$, where the step of $\varphi_l$ equals $l$.
The notion of the step is well-defined on $\gt z(e)$, due to
an equality $(-i)-(-j)=j-i$. From the definition of
$\alpha$, one can deduce that
$\alpha(\varphi_l,\varphi_t)\ne 0$ only if $l+t=1$. The
stabiliser $\gt z(e)_\alpha$ is a direct sum of its
subspaces $\Phi_{l}$, consisting of elements having step
$l$. As we have seen, $\Phi_l=\varnothing$ if $l>0$. It is
remains to describe elements with non-positive steps.

\begin{ex} Let us show that $\dim\Phi_0\le d_0/2$.
Suppose  $\varphi\in\Phi_0$, $\varphi\ne 0$ and
$\varphi(w_0)=0$. Take a minimal by the absolute value $i$
such that $\varphi(w_i)\ne 0$. Since $\varphi\in\gt{so}_n$,
we have also $\varphi(w_{-i})\ne 0$. Assume that $i>0$ and a
coefficient $c_i^{i,s}$ of $\varphi$ is non-zero. Then
$|i-1|<i$, $d_{i-1}\ge d_i$, there is an element
$\xi_{i-1}^{i,d_i-s}\in\gt z_{\gt{gl}}(e)$ and
$0=\alpha([\xi,\varphi])=c_i^{i,s}-c_{i-1}^{i-1,s}=c_i^{i,s}$.
Hence, if $\varphi(w_0)=0$, then also $\varphi=0$. Thus, a
vector $\varphi\in\Phi_0$ is entirely determined by its
value on $w_0$. In its turn
$\varphi(w_0)=c_1ew_0+c_3e^3w_0+\ldots+
c_{d_0-1}e^{d_0-1}w_0$.
\end{ex}

\begin{lm} If $q=2l$ or $q=2l-1$, where $0<l\le m$, then
$\dim\Phi_{-q}\le (d_l+1)/2$.
\end{lm}

\begin{proof}
Similarly to the previous example, we show that if
$\varphi\in\Phi_{-q}$ and $\varphi(w_l)=0$, then also
$\varphi=0$. Since $\varphi\in\gt{so}_n$, if
$\varphi(w_i)\ne 0$, then also $\varphi(w_{q-i})\ne 0$.
Suppose $\varphi(w_j)\ne 0$ for some $j$. If $j<l$, then
$j-q\ge l$, but $\varphi(w_l)=0$, hence $j>l$. Find the
minimal $j>l$ such that $\varphi(w_j)\ne 0$. Suppose
$c_j^{j-q,s}=c_j^{j-q,s}(\varphi)\ne 0$. We have $-j<-l\le
j-q-1<j$, $d_j\le d_{j-q-1}$, $d_{j-q}-s\ge 0$. Hence, there
is an element $\xi:=\xi_{j-q-1}^{j,d_{j-q}-s}\in\gt
z_{\gt{gl}}(e)$. As above
$0=\alpha([\xi,\varphi,])=c_j^{j-q,s}-
c_{j-1}^{j-q-1,\delta}=c_j^{j-q,s}$ (we do not give a
precise value of $\delta$, anyway, $\varphi(w_{j-1})=0$,
since $l\le j-1<j$). To conclude we describe possible values
$\varphi(w_l)$. If  $q=2l$, then
$\varphi(w_l)=c_0w_{-l}+c_2e^2w_{-l}+\ldots+
c_{d_l}e^{d_l}w_{-l}$. In case $q=2l-1$ we get an equation
on coefficients of $\varphi$:
$0=\alpha[(\xi_{-l}^{l,b},\varphi])=
c_l^{-l+1,d_{l-1}-b}-c_{l-1}^{-l,d_l-b}$, i.e.,
$c_l^{-l+1,d_{l-1}-b}=c_{l-1}^{-l,d_l-b}$. This is possible
only for odd $b$.
\end{proof}

\begin{thm} Suppose $e\in\gt{so}_n$ is a nilpotent element.
Then $\ind\gt z(e)=\rk\gt{so}_n=[n/2]$.
\end{thm}

\begin{proof} If possibility
{\bf (3)} takes place, i.e., only one Jordan block of $e$ is
odd-dimensional and it is also maximal, then, as we have
seen, $\gt
z(e)_{\tilde\alpha}=\bigoplus\limits_{q=0}^{2m}\Phi_{-q}$.
Moreover, $\dim\Phi_q$ is at most half of dimension of the
Jordan block with number $[(q+1)/2]$. Thereby, $\dim\gt
z(e)_{\tilde\alpha}\le((d_0+1)/2+(\sum\limits_{l=1}^{m}
d_l))=[n/2]$. On the other hand, according to Vinberg's
inequality, $\ind\gt z(e)\ge [n/2]$.

In cases {\bf (1)} and {\bf (2)} the inequality $\ind\gt
z(e)\le\rk\gt{so}_n$ was proved by induction.

If none of these three possibilities takes place, then
either $k=1$ and $e$ is a regular nilpotent element, or
$k=2$ and both Jordan blocks of $e$ are even-dimensional.
These two  cases  have been considered separately.
\end{proof}

\section{Generic points}
In this section we assume that
$\mathbb K$ is algebraically closed
and of characteristic zero. Suppose we have a linear
action of a Lie algebra
$\gt g$ on a vector space $V$.

\vskip0.2ex {\bf Definition.} A vector $x\in V$ (a
subalgebra $\gt g_x$) is called {\it a generic point} ({\it
a generic stabiliser}), if for every point $y\in U\subset V$
of some open in Zariski topology subset $U$ algebras $\gt
g_y$ and $\gt g_x$ are conjugated in $\gt g$. \vskip0.4ex

\noindent
It is well known that generic points exist
for any linear action of a reductive Lie algebra.

It is proved in \cite[\textsection 1]{al2} that
a subalgebra
$\gt g_x$ is a generic stabiliser
\iff $V=V^{\gt g_x}+\gt gx$, where
$V^{\gt g_x}$ is the subspace of all
vectors of $V$
invariant under $\gt g_x$.

Tauvel and Yu have noticed that in case of a coadjoint
representation $\gt gx=(\gt g/\gt g_x)^*=\Ann(\gt g_x)$,
$(\gt g^*)^{\gt g_x}=\Ann([\gt g_x, \gt g])$. From this
observation they have deduced a simple and useful criterion.

\begin{thm}\label{yu}\cite[Corollaire 1.8.]{yu}
Let $\gt g$ be a Lie algebra and $x\in\gt g^*$. The subalgebra $\gt
g_x$ is a generic stabiliser of the coadjoint
representation of $\gt g$ \iff
$[\gt g_x,\gt g]\cap\gt g_x=\{0\}$.
\end{thm}

\noindent Unfortunately, the authors of \cite{yu} were not
aware of the aforementioned \'Elashvili's result and have
proved it anew.

Let $e\in\gt{gl}_n$ be a nilpotent element
and
$\gt z(e)$ the centraliser of $e$. Set
 $\gt h=\gt z(e)_\alpha$, where
$\alpha\in\gt z(e)^*$ is the same as in Section~2.

\begin{prop} There is an
$\gt h$-invariant decomposition
$\gt z(e)=\gt h\oplus\gt m$, where
$\gt m$ is generated by the
vectors $\xi_i^{j,s}$ with $i\ne j$.
\end{prop}

\begin{proof} Recall that $\gt h$ is generated by
the vectors  $\xi_i^{i,s}$.
The inclusion $[\gt h,\gt m]\subset\gt m$
follows immediately from the equality
$$
[\xi_i^{i,s},\xi_j^{t,b}]=\left\{
\begin{array}{rl}
\xi_i^{t,s+b} &\mbox{ if } i=j, i\ne t; \\
-\xi_j^{i,s+b} &\mbox{ if } i=t, i\ne j;\\
0 &\mbox{ otherwise.} \\
\end{array}\right.
$$
\end{proof}

There is a similar decomposition in the case of symplectic
algebras. Let $e\in\gt{sp}(V)\subset\gt{gl}(V)$. Denote by
$\gt z_{\gt{gl}}(e)$ and $\gt z_{\gt{sp}}(e)$ the
centralisers of $e$ in $\gt{gl}(V)$ and $\gt{sp}(V)$,
respectively. We use notation of Section~3. Suppose $\gt
z_{\gt{gl}}(e)=\gt h\oplus\gt m$. Evidently, this
decomposition is $\sigma$-invariant and $\gt
z_{\gt{sp}}(e)=\gt h^\sigma\oplus\gt m^\sigma$, where $\gt
h^\sigma=\gt z_{\gt{sp}}(e)_{\tilde\alpha}$.

\begin{thm} The Lie algebras $\gt
z_{\gt{gl}}(e)_\alpha$ and $\gt
z_{\gt{sp}}(e)_{\tilde\alpha}$  constructed in Sections~2
and 3 in cases of general linear and symplectic algebras are
generic stabilisers of the coadjoint actions of $\gt
z_{\gt{gl}}(e)$ and $\gt z_{\gt{sp}}(e)$.
\end{thm}

\begin{proof} Let us verify the condition of Theorem~\ref{yu}.
Since
$[\gt h, \gt z_{\gt{gl}}(e)]=[\gt h,\gt m]\subset\gt m$, we have
 $[\gt h, \gt z_{\gt{gl}}(e)]\cap\gt h=0$. Similarly,
$[\gt h^\sigma, \gt z_{\gt{sp}}(e)]\subset\gt m^\sigma$.
\end{proof}

In case of orthogonal algebras it can happen that
a generic stabiliser
of the coadjoint action of
$\gt z(e)$ does not exist.

\begin{ex} Let $e\in\gt{so}_8$ be a subregular nilpotent
element. Then it has two Jordan blocks of dimensions
3 and 5.
Choose the generators $w_1$, $w_2$ of Jordan blocks such that
$(w_1, e^2w_1)_V=(w_2,e^4w_2)_V=1$. The dimension of $\gt
z(e)$ is 6 and $\gt z(e)$ has a three-dimensional center,
generated by the vectors $e$, $e^3=\xi_2^{2,3}$ and
$\varphi_3=\xi_1^{2,4}-\xi_2^{1,2}$. Since $\ind\gt z(e)=4$,
we have $\dim\gt z(e)_\alpha=4$ for points of some open
subset $U\subset\gt z(e)^*$.

Assume that a generic stabiliser of the coadjoint action of
$\gt z(e)$ exists and denote it by $\gt f$. Evidently, $\gt
f$ contains the center of  $\gt z(e)$. Consider an element
$\varphi_2=\xi_1^{2,3}+\xi_2^{1,1}\in\gt z(e)$. Clearly, the
subspace $[\varphi_2, \gt z(e)]$ is a linear span of $e^3$
and $\varphi_3$. In particular, it is contained in the
center of $\gt z(e)$, and, hence, in $\gt f$. Hence,
$[\varphi_2, \gt f]\subset\gt f$, and, by Theorem~\ref{yu},
$\gt f\subset\gt z(e)_{\varphi_2}$. Since $\dim\gt
z(e)_{\varphi_2}=4$, we have $\gt f=\gt z(e)_{\varphi_2}$.
On the other hand, $\gt z(e)_{\varphi_2}=\left<e,e^3,
\varphi_3, \varphi_2\right>_{\mathbb K}$ is a normal, but
not a central subalgebra of $\gt z(e)$.

Consider the embedding
$\gt{so}_8\subset\gt{so}_9$ as the stabiliser
of the first basis vector in  $\mathbb K^9$.
By a similar argument
one can show that a generic stabiliser does not exist
for the coadjoint action of
$\gt z_{\gt{so}_9}(e)$ either.
\end{ex}

\section{Commuting varieties}

Let $\gt g$ be a Lie algebra over an algebraically closed
field $\mathbb K$ of characteristic zero. A closed subset
$Y=\{(x,y)|x,y\in\gt g, [x,y]=0\}\subset(\gt g\times\gt g)$
is called the {\it commuting variety} of the algebra $\gt
g$. The question of whether $Y$ is irreducible or not is of a
great interest. In case of a reductive algebra $\gt g$ the
commuting variety  $Y$ is irreducible and coincides with
the closure of
${G(\gt a, \gt a)}$, where $\gt a\subset\gt
g$ is a Cartan subalgebra and $G$ is a connected algebraic group with
$\Lie G=\gt g$.

Let $e\in\gt{gl}_n$ be a nilpotent element and $\gt z(e)$
the centraliser of $e$. We use notation introduced in
Section 2. Set $\gt h=\gt z(e)_\alpha$. Consider a
subalgebra $\gt t\subset\gt z(e)$ generated by the vectors
$\xi_i^{i,0}$. Evidently, $\gt t\subset\gt h$. Moreover,
since $[\xi_i^{j,s}, t_i\xi_i^{i,0}+t_j\xi_j^{j,0}]=
(t_j-t_i)\xi_i^{j,s}$, the algebra $\gt h$ coincides with
the normaliser (= centraliser) of $\gt t$ in $\gt z(e)$.
Hence, $\gt h$ coincides with its normaliser in $\gt z(e)$.

Let $Z(e)$ be the identity component of the centraliser of
$e$ in ${\rm GL}_n$. Then $Y_0=\overline{Z(e)(\gt h, \gt
h)}$ is an irreducible component of $Y$ of maximal
dimension. As in the reductive case, $Y$ is irreducible \iff
$Y_0=Y$. It is known that if a nilpotent element $e$ has at
most two Jordan blocks, then $Y$ is irreducible \cite{2}. In
the general case, the statement is not true, since it would
lead to the irreducibility of the commuting varieties of
triples of matrices.

\begin{ex} Assume that
$Y_0=Y$ for all nilpotent elements $e\in\gt{gl}_m$ with
$m\le n$. Consider the set of triples of commuting matrices
$$
C_3=\{(A,B,C)|A,B,C\in\gt{gl}_n, [A,B]=[A,C]=[B,C]=0\}.
$$
Let
$\gt a\subset\gt{gl}_n$ be a subalgebra of diagonal
matrices. Clearly,
$\overline{\GL_n(\gt a, \gt a, \gt a)}$ is an irreducible component
of $C_3$. Let us prove by induction that it coincides with
$C_3$. There is nothing to prove for
$n=1$. Let $n>1$.
We show that each triple $(A,B,C)$ of commuting
matrices is contained in the closure
$\overline{\GL_n(\gt a, \gt a, \gt a)}$.
Without loss of generality,  we may assume that
$A,B,C\in\gt{sl}_n$. Let $A=A_s+A_n$ be the Jordan decomposition
of $A$. If $A_s\ne 0$, consider the centraliser $\gt z(A_s)$ of $A_s$
in $\gt{gl}_n$. Clearly,
$A,B,C\in\gt z(A_s)$ and $\gt z(A_s)$ is a sum
of several algebras
$\gt{gl}_{n_i}$ with strictly smaller dimension.
We may assume that
$\gt a\subset\gt z(A_s)$. Then, by
the inductive hypothesis
$$
(A,B,C)\in\overline{Z(A_s)(\gt a,\gt a,\gt a)}\subset
\overline{\GL_n(\gt a, \gt a, \gt a)}.
$$

Suppose now that all three elements $A,B,C$ are nilpotent
and at least one of them, say $A$, is not regular.
Consider the centraliser $\gt z(A)\subset\gt{gl}_n$. We have
assumed that $Y_0=Y$, i.e., the pair $(B,C)$ lies in the
closure of 
$Z(A)(\gt h,\gt h)$. It will be enough to show that $(A,\gt
h,\gt h)\subset\overline{\GL_n(\gt a, \gt a, \gt a)}$. Let
$x\in\gt t\subset\gt h$ be a non-central semisimple element.
Then $A\in(\gt{gl}_n)_x$ and $\gt h\subset(\gt{gl}_n)_x$.
Once again we can make an induction step, passing to a
subalgebra $(\gt{gl}_n)_x$.

If all three elements $A,B,C$ are regular nilpotent, then
there is a non-trivial linear combination $A'$ of them,
which is non-regular. In particular, the triple $(A,B,C)$ is
equivalent under the action of $\GL_n$ to some other triple
$(A',B',C')$ of commuting nilpotent matrices.
\end{ex}

It is known that for $n>31$ the variety $C_3$ is reducible,
see  \cite{gur}. Hence,  
the commuting variety $Y$
is certainly reducible for some nilpotent elements. It will be
interesting to find minimal (in some sense) nilpotent
elements for which $Y$ is reducible and/or describe some
classes of nilpotent elements for which
 $Y$ is irreducible.

\end{document}